\documentclass[11pt]{article}
\usepackage{times}
\usepackage{indentfirst}
\usepackage{amsthm,amsmath}
\usepackage{calrsfs}

\newtheorem{Theorem}{Theorem}
\newtheorem{Lemma}{Lemma}
\newtheorem{Proposition}{Proposition}

\newcounter{remark}

\newcounter{definition}

\newcounter{example}

\begin{document}

\section*{\centering\Large\rm
          Lower limits for distributions of randomly stopped sums
%\\
 %         in the presence of heavy tails
          \footnote{The research of Denisov and Foss
          was partially supported by the EPSRC  Grant
          EP/E033717/1.
          The research of Foss and Korshunov
          was partially supported by the Royal Society
          International Joint Project Grant 2005/R2-JP.
          The research of Korshunov was partially
          supported by Russian Science Support Foundation.}}
\section*{\centering\large\rm
D. Denisov,\footnote{Address: School of MACS,
Heriot-Watt University,
Edinburgh EH14 4AS, UK.
E-mail address: Denisov@ma.hw.ac.uk}
S. Foss,\footnote{Address: School of MACS,
Heriot-Watt University,
Edinburgh EH14 4AS, UK;
and Sobolev Institute of Mathematics,
4 Koptyuga Pr., Novosibirsk 630090, Russia.
E-mail address: S.Foss@ma.hw.ac.uk}
and D. Korshunov\footnote{Address:
Sobolev Institute of Mathematics,
4 Koptyuga pr., Novosibirsk 630090, Russia.
E-mail address: Korshunov@math.nsc.ru}}
\section*{\centering\normalsize\it
Eurandom, Heriot-Watt University, and Sobolev Institute of Mathematics}

\begin{abstract}
We study lower limits for the ratio
$\frac{\overline{F^{*\tau}}(x)}{\overline F(x)}$
of tail distributions where $ F^{*\tau}$ is a distribution
of a sum of a random size $\tau$ of
i.i.d. random variables having a common distribution $F$,
and a random variable $\tau$ does not depend on summands.

AMS classification: Primary 60E05; secondary 60F10

{\it Keywords:}
Convolution tail; Random sums of random variables;
Lower limit; Heavy- and light-tailed distributions
\end{abstract}

\hspace{20mm}

{\bf 1. Introduction.}
Let $\xi$, $\xi_1$, $\xi_2$, \ldots~be independent identically
distributed random variables.
We assume that their common distribution $F$
is unbounded from the right, that is,
$\overline F(x)\equiv F(x,\infty)>0$ for all $x$.
Put $S_0=0$ and $S_n=\xi_1+\ldots+\xi_n$, $n=1$, $2$, \ldots.

Let $\tau$ be a counting random variable which
does not depend on $\{\xi_n\}_{n\ge 1}$.
Denote by $F^{*\tau}$ the distribution of a random
sum $S_\tau=\xi_1+\ldots+\xi_\tau$.
In this paper we study lower limits (as $x\to\infty$)
for the ratio
$\frac{\overline{F^{*\tau}}(x)}{\overline F(x)}$.

We distinguish two types of distributions, heavy- and light-tailed.
A random variable $\eta$ has a {\it heavy-tailed} distribution if
${\bf E} e^{\varepsilon\eta} = \infty$ for all
$\varepsilon>0$, and {\it light-tailed} otherwise.

We consider only non-negative random variables
and, in the case of heavy-tailed $F$,
study conditions for
\begin{eqnarray}\label{lower.limit.concl}
\liminf_{x\to\infty}\frac{\overline{F^{*\tau}}(x)}{\overline F(x)}
&=& {\mathbf E}\tau
\end{eqnarray}
to hold.
This problem has been given a complete solution
in [\ref{FK}] for $\tau=2$, and then
in [\ref{DFK2}] for  $\tau$  with a light-tailed
distribution and for heavy-tailed summands.
In the present work, we generalise results of [\ref{DFK2}] onto
classes
of distributions of $\tau$ which include all light-tailed
distributions
and also some heavy-tailed distributions.
With each heavy-tailed distribution $F$,
we associate a corresponding class of distributions of $\tau$.
For earlier studies on lower limits and
on a related problem of justifying a constant $K$
in the equivalence
$\overline{F^{*2}}(x) \sim K \overline F(x)$, see e.g.
[\ref{CNW}, \ref{Cline}, \ref{EG}, \ref{Pakes}, \ref{Rogozin}]
and further references therein.

Since the inequality ``$\ge$'' in (\ref{lower.limit.concl}) is valid
for
non-negative
$\{\xi_n\}$ without any further assumptions
(see, e.g., [\ref{Rudin}] or [\ref{DFK2}]),
we immediately get the equality if ${\bf E}\tau=\infty$.
Therefore, in the rest of the paper, we consider the case
${\bf E}\tau<\infty$ only. Our first result is

\begin{Theorem}\label{ht.th}
Assume that $\xi\ge0$ is heavy-tailed
and ${\bf E}\xi<\infty$.
Let, for some $c>{\bf E}\xi$,
\begin{eqnarray}\label{G.=o.F}
{\bf P}\{c\tau>x\} &=& o(\overline F(x))
\quad\mbox{ as }x\to\infty.
\end{eqnarray}
Then {\rm(\ref{lower.limit.concl})} holds.
\end{Theorem}

The proof of Theorem 1 is based on a study
of moments ${\bf E}e^{f(\xi)}$
for appropriately chosen concave function $f$.
More precisely, we deduce Theorem \ref{ht.th}
from the following general result which
explores some ideas from
[\ref{Rudin}, \ref{FK}, \ref{DFK2}].

\begin{Theorem}\label{main_result}
Assume that $\xi\ge0$ is heavy-tailed
and ${\bf E}\xi<\infty$.
Let there exists a function
$f:{\bf R}^+\to{\bf R}$ such that
\begin{eqnarray}\label{eq.1}
{\bf E}e^{f(\xi)}=\infty,
\end{eqnarray}
and, for some $c>{\bf E}\xi$,
\begin{eqnarray}\label{eq.2}
{\bf E}e^{f(c\tau)}<\infty.
\end{eqnarray}
If $f(x)\ge\ln x$ for all sufficiently large $x$ and if
the difference $f(x)-\ln x$ is an eventually concave function,
then {\rm(\ref{lower.limit.concl})} holds.
\end{Theorem}

In particular, the equality (\ref{lower.limit.concl}) is valid
provided ${\bf E}\xi^k=\infty$ and ${\bf E}\tau^k<\infty$
for some $k\ge1$; it is sufficient to consider
the function $f(x)=k\ln x$. Earlier this was proved in
[\ref{DFK2}, Theorem 1] by a more simple method.

If we consider instead the function $f(x)=\gamma x$, $\gamma>0$,
then we obtain the equality (\ref{lower.limit.concl})
provided $\xi$ is heavy-tailed but $\tau$ is light-tailed.
This is Theorem 2 from [\ref{DFK2}].

Finally, the equality (\ref{lower.limit.concl}) is valid
if $F$ is a Weibull distribution with parameter $\beta \in (0,1)$,
$\overline F(x)=e^{-x^\beta}$ and $f(x)=x^\beta$ or,
more generally, $f(x)=x^\beta - c\ln x$ for $x\ge1$
where $c\le\beta$ is any fixed constant.

The counterpart of Theorem \ref{ht.th}
in the light-tailed case is stated next.
But first we need some notations.
By the Laplace transform of $F$ at the point
$\gamma\in{\bf R}$ we mean
$$
\varphi(\gamma)=\int_0^\infty e^{\gamma x}F(dx)
\in(0,\infty].
$$
Put
$$
\widehat\gamma=\sup\{\gamma:\varphi(\gamma)<\infty\}
\in[0,\infty].
$$
Note that the function $\varphi(\gamma)$ is monotone
continuous in the interval $(-\infty,\widehat\gamma)$,
and $\varphi(\widehat\gamma)=
\lim\limits_{\gamma\uparrow\widehat\gamma}\varphi(\gamma)\in[1,\infty]$.

\begin{Theorem}\label{lt.th}
Let $\widehat\gamma\in(0,\infty]$, so that
$\varphi(\widehat\gamma)\in(1,\infty]$.
If {\rm(\ref{G.=o.F})} holds and, for any fixed $y>0$,
\begin{eqnarray}\label{R.R.1.th}
\liminf_{x\to\infty}
\frac{\overline F(x-y)}{\overline F(x)}
&\ge& e^{\widehat\gamma y},
\end{eqnarray}
then
\begin{eqnarray*}
\liminf_{x\to\infty}
\frac{\overline{F^{*\tau}}(x)}{\overline F(x)}
&=& {\bf E}\tau\varphi^{\tau-1}(\widehat\gamma).
\end{eqnarray*}
\end{Theorem}

The paper is organised as follows. In Section 2,
we formulate and prove a general result
on characterisation of heavy-tailed distributions
on the positive half-line.
Section 3 is devoted to the estimation
of the functional ${\bf E}e^{h(S_n)}$
for a concave function $h$.
Sections 4 and 5 contain proofs of
Theorems \ref{main_result} and \ref{ht.th}
respectively.
Section 6 is devoted to the proof in light-tailed case.

{\bf 2. Characterisation of heavy-tailed distributions.}
It was proved in [\ref{DFK2}, Lemma 2] that,
for any heavy-tailed random variable $\xi\ge 0$
and for any real $\delta>0$,
there exists an increasing concave
function $h:{\bf R}^+\to{\bf R}^+$ such that
${\bf E}e^{h(\xi)}\le 1+\delta$ and
${\bf E}\xi e^{h(\xi)}=\infty$.
In the present section, we obtain some generalisation of it.

\begin{Lemma}\label{lm.moments.ext}
Let $\xi\ge 0$ be a random variable
with a heavy-tailed distribution.
Let $f:{\bf R}^+\to{\bf R}$
be a concave function such that
\begin{eqnarray}\label{inf.f}
{\bf E}e^{f(\xi)} = \infty.
\end{eqnarray}
Let a function $g:{\bf R}^+\to{\bf R}$
be such that $g(x)\to\infty$ as $x\to\infty$.
Then there exists a concave function
$h:{\bf R}^+\to{\bf R}$ such that $h\le f$ and
\begin{eqnarray*}
{\bf E}e^{h(\xi)}<\infty,
&& {\bf E}e^{h(\xi)+g(\xi)}=\infty.
\end{eqnarray*}
\end{Lemma}

\proof
Without loss of generality assume $f(0)=0$.
We will construct a function $h(x)$ on the
successive intervals.
For that we introduce two positive sequences,
$x_n\uparrow\infty$ as $n\to\infty$
and $\varepsilon_n\in(0,1]$.
We put $x_0=0$, $h(0)=f(0)=0$, $h'(0)=f'(0)$, and
$$
h(x)=h(x_{n-1})+\varepsilon_n
\min(h'(x_{n-1})(x-x_{n-1}),f(x)-f(x_{n-1}))
\quad\mbox{for }x\in(x_{n-1},x_n];
$$
here $h'$ is the left derivative of the function $h$.
The function $h$ is increasing, since $\varepsilon_n>0$
and $f$ is increasing.
Moreover, this function is concave, due to
$\varepsilon_n\le 1$ and concavity of $f$.
Since $h(x)-h(x_{n-1})\le f(x)-f(x_{n-1})$
for $x\in(x_{n-1},x_n]$, we have $h\le f$.

Now proceed with the very construction
of $x_n$ and $\varepsilon_n$.
By conditions $g(x)\to\infty$ and (\ref{inf.f}),
we can choose $x_1$ so large that
$e^{g(x)}\ge 2^1$ for all $x\ge x_1$ and
\begin{eqnarray*}
{\bf E}\{e^{\min(h'(0)\xi,f(\xi))};\xi\in(x_0,x_1]\}
+e^{\min(h'(0)x_1,f(x_1))}\overline F(x_1)
&>& \overline F(x_0)+1.
\end{eqnarray*}
Choose $\varepsilon_1\in(0,1]$ so that
\begin{eqnarray*}
{\bf E}\{e^{\varepsilon_1\min(h'(0)\xi,f(\xi))}; \xi\in(x_0,x_1]\}
+e^{\varepsilon_1\min(h'(0)x_1,f(x_1))}\overline F(x_1)
&=& \overline F(x_0)+1.
\end{eqnarray*}
Put $h(x)=\varepsilon_1\min(x, f(x))$ for $x\in(0,x_1]$.
Then the latter equality is equivalent to
\begin{eqnarray*}
{\bf E}\{e^{h(\xi)}; \xi\in(x_0,x_1]\}
+e^{h(x_1)}\overline F(x_1)
&=& e^{h(x_0)}\overline F(x_0)+1/2,
\end{eqnarray*}

By induction we construct an increasing sequence
$x_n$ and a sequence $\varepsilon_n\in(0,1]$
such that $e^{g(x)}\ge 2^n$ for all $x\ge x_n$, and
\begin{eqnarray*}
{\bf E}\{e^{h(\xi)}; \xi\in(x_{n-1},x_n]\}
+e^{h(x_n)}\overline F(x_n)
&=& e^{h(x_{n-1})}\overline F(x_{n-1})+1/2^n
\end{eqnarray*}
for any $n\ge 1$.
For $n=1$ this is already done.
Make the induction hypothesis for some $n\ge2$.
For any $x>x_n$, denote
\begin{eqnarray*}
\delta(x,\varepsilon) &\equiv&
e^{h(x_n)}\Bigl({\bf E}\{e^{\varepsilon\min(h'(x_n)(\xi-x_n),
f(\xi)-f(x_n))}; \xi\in(x_n,x]\}\\
&& \hspace{60mm} +e^{\varepsilon\min(h'(x_n)(x-x_n),f(x)-f(x_n))}
\overline F(x)\Bigr).
\end{eqnarray*}
By the convergence $g(x)\to\infty$,
by heavy-tailedness of $\xi$,
and by the condition (\ref{inf.f}),
there exists $x_{n+1}$ so large that
$e^{g(x)}\ge 2^{n+1}$ for all $x\ge x_{n+1}$ and
\begin{eqnarray*}
\delta(x_{n+1},1)
&>& e^{h(x_n)}\overline F(x_n)+1.
\end{eqnarray*}
Note that the function
$\delta(x_{n+1},\varepsilon)$
is continuously decreasing to $e^{h(x_n)}\overline F(x_n)$
as $\varepsilon\downarrow 0$.
Therefore, we can choose
$\varepsilon_{n+1}\in(0,1]$ so that
\begin{eqnarray*}
\delta(x_{n+1},\varepsilon_{n+1})
&=& e^{h(x_n)}\overline F(x_n)+1/2^{n+1}.
\end{eqnarray*}
Then
\begin{eqnarray*}
{\bf E}\{e^{h(\xi)};\xi\in(x_n,x_{n+1}]\}
+e^{h(x_{n+1})}\overline F(x_{n+1})
&=& e^{h(x_n)}\overline F(x_n)+1/2^{n+1}.
\end{eqnarray*}
Our induction hypothesis now holds with $n+1$
in place of $n$ as required.

Next, for any $N$,
\begin{eqnarray*}
{\bf E}\{e^{h(\xi)};\xi\le x_{N+1}\}
&=& \sum_{n=0}^N
{\bf E}\{e^{h(\xi)}; \xi\in(x_n,x_{n+1}]\}\\
&=& \sum_{n=0}^N
\Bigl(e^{h(x_n)}\overline F(x_n)-
e^{h(x_{n+1})}\overline F(x_{n+1})+1/2^{n+1}\Bigr)\\
&\le& e^{h(x_0)}\overline F(x_0)+1,
\end{eqnarray*}
so that ${\bf E}e^{h(\xi)}$ is finite.
On the other hand, since
$e^{g(x)}\ge 2^k$ for all $x\ge x_k$,
\begin{eqnarray*}
{\bf E}\{e^{h(\xi)+g(\xi)};\xi>x_n\}
&\ge& 2^n\Bigl({\bf E}\{e^{h(\xi)}; \xi\in(x_n,x_{n+1}]\}
+e^{h(x_{n+1})}\overline F(x_{n+1})\Bigr)\\
&=& 2^n(e^{h(x_n)}\overline F(x_n)+1/2^{n+1}).
\end{eqnarray*}
Then, for any $n$,
${\bf E}\{e^{h(\xi)+g(\xi)};\xi>x_n\} \ge 1/2$,
which implies ${\bf E}e^{h(\xi)+g(\xi)}=\infty$.
The proof is complete.

\begin{Lemma}\label{lm.moments.ext.f1f2}
Let $\xi\ge 0$ be a random variable with a heavy-tailed distribution.
Let $f_1:{\bf R}^+\to{\bf R}$ be any measurable
function  and $f_2:{\bf R}^+\to{\bf R}$
a concave function such that
\begin{eqnarray*}
{\bf E}e^{f_1(\xi)} < \infty
&\mbox{ and }& {\bf E}e^{f_1(\xi)+f_2(\xi)} = \infty.
\end{eqnarray*}
Let a function $g:{\bf R}^+\to{\bf R}$
be such that $g(x)\to\infty$ as $x\to\infty$.
Then there exists a concave function
$h:{\bf R}^+\to{\bf R}$ such that $h\le f_2$ and
\begin{eqnarray*}
{\bf E}e^{f_1(\xi)+h(\xi)}<\infty
&\mbox{ and }& {\bf E}e^{f_1(\xi)+h(\xi)+g(\xi)}=\infty.
\end{eqnarray*}
\end{Lemma}

\proof
Consider a new governing probability measure
${\bf P}^*$ defined in the following way:
$$
{\bf P}^*\{d\omega\}
= \frac{e^{f_1(\xi(\omega))}{\bf P}\{d\omega\}}
{{\bf E}e^{f_1(\xi)}}.
$$
Then
$$
{\bf E}^*e^{f_2(\xi)}
=\frac{{\bf E}e^{f_1(\xi)+f_2(\xi)}}
{{\bf E}e^{f_1(\xi)}}
=\infty.
$$
In particular, $\xi$ is heavy-tailed against
the measure ${\bf P}^*$.
Now it follows from Lemma \ref{lm.moments.ext}
that there exists a concave function
$h:{\bf R}^+\to{\bf R}$ such that $h\le f_2$,
$h(x)=o(x)$,
${\bf E}^*e^{h(\xi)}<\infty$, and
${\bf E}^*e^{h(\xi)+g(\xi)}=\infty$.
Equivalently,
$$
{\bf E}e^{f_1(\xi)+h(\xi)}
={\bf E}e^{f_1(\xi)} {\bf E}^*e^{h(\xi)}
<\infty
$$
and
$$
{\bf E}e^{f_1(\xi)+h(\xi)+g(\xi)}
={\bf E}e^{f_1(\xi)} {\bf E}^*e^{h(\xi)+g(\xi)}
=\infty.
$$
The proof is complete.

{\bf 3. Growth rate of sums in terms of generalised moments.}
According to the Law of Large Numbers,
the sum $S_n$ growths like $n{\bf E}\xi$.
In the following lemma we provide conditions
on a function $h(x)$, guaranteeing an appropriate
rate of growth for the functional ${\bf E}e^{h(S_n)}$.

\begin{Lemma}\label{lem.conc2}
Let $\xi$ be a non-negative random variable.
Let $h:{\bf R}^+\to{\bf R}$ be a non-decreasing
eventually concave function such that
$h(x)=o(x)$ as $x\to\infty$ and $h(x)\ge\ln x$
for all sufficiently large $x$.
If ${\bf E}e^{h(\xi)}<\infty$, then, for any
$c>{\bf E}\xi$, there exists a constant $K(c)$
such that ${\bf E}e^{h(S_n)} \le K(c) e^{h(nc)}$,
for all $n$.
\end{Lemma}

To prove this lemma, we need the following assertion,
which generalises the corresponding estimate from [\ref{FS}]:

\begin{Lemma}\label{semi.moments.th}
Let $\eta$ be a random variable with ${\bf E}\eta<0$.
Let $h:{\bf R}\to{\bf R}$ be a non-decreasing
and eventually concave function such that
$h(x)=o(x)$ as $x\to\infty$ and $h(x)\ge\ln x$
for all sufficiently large $x$.
If ${\bf E}e^{h(\eta)}<\infty$,
then there exists $x_0$ such that the inequality
${\bf E}e^{h(x+\eta)}\le e^{h(x)}$
holds for all $x>x_0$.
\end{Lemma}

\proof
Since $h$ is increasing, without loss of
generality we may assume that $\eta$ is bounded
from below, that is, $\eta\ge M$ for some $M$.
Also, we may assume that $h$ is non-negative and
concave on the whole half-line $[0,\infty)$.

Since $h$ is concave, $h'(x)$ is non-increasing function.
With necessity $h'(x)\to 0$ as $x\to\infty$,
otherwise the condition $h(x)=o(x)$ is violated.
If ultimately $h'(x)=0$, then $h$ is
ultimately a constant function and
the proof of the theorem is obvious.

Consider now the case $h'(x)\to 0$
as $x\to\infty$ but $h'(x)>0$ for all $x$.
Put $g(x)\equiv 1/h'(x)$,
then $g(x)\uparrow\infty$ as $x\to\infty$.
Since ${\bf E}\eta<0$, we can choose
sufficiently large $A$ such that
\begin{eqnarray}\label{choice.A}
\varepsilon \equiv
{\bf E}\{\eta;\eta\in[M,A]\}
+e{\bf E}\{\eta;\eta>A\} &<& 0.
\end{eqnarray}
By concavity of $h$, for any $x$ and
$y\in{\bf R}$ we have the inequality
$h(x+y)-h(x) \le h'(x)y$. Hence,
\begin{eqnarray}\label{decomp.estimate}
{\bf E}e^{h(x+\eta)-h(x)}
&\le& {\bf E}\{e^{h'(x)\eta}; \eta\in[M,A]\}
+{\bf E}\{e^{h'(x)\eta}; \eta\in(A,g(x)]\}
\nonumber\\
&& +{\bf E}\{e^{h(x+\eta)-h(x)};\eta>g(x)\}\nonumber\\
&\equiv& E_1+E_2+E_3.
\end{eqnarray}
Since $h'(x)\to 0$,
the Taylor's expansion for the exponent
up to the linear term implies, as $x\to\infty$,
\begin{eqnarray}\label{E1}
E_1 &=&
{\bf P}\{\eta\in[M,A]\}
+h'(x){\bf E}\{\eta;\eta\in[M,A]\}+o(h'(x)).
\end{eqnarray}
On the event $\eta\in(A,g(x)]$
we have $h'(x)\eta\le 1$ and, thus,
$e^{h'(x)\eta}\le 1+eh'(x)\eta$. Then
\begin{eqnarray}\label{E2}
E_2 &\le& {\bf P}\{\eta\in(A,g(x)]\}
+eh'(x){\bf E}\{\eta;\eta\in(A,g(x)]\}.
\end{eqnarray}
We have
\begin{eqnarray}\label{E3.1}
E_3 &=& {\bf E}\{e^{h(\eta)}
e^{h(x+\eta)-h(x)-h(\eta)};\eta>g(x)\}.
\end{eqnarray}
By concavity of $h$, for $x>0$, the difference
$h(x+y)-h(y)$ is non-increasing in $y$. Therefore,
for any $y>g(x)$,
\begin{eqnarray*}
h(x+y)-h(x)-h(y)
&\le& h(x+g(x))-h(x)-h(g(x))\\
&\le& h'(x)g(x)-h(g(x))\\
&=& 1-h(g(x))\\
&\le& 1-\ln g(x),
\end{eqnarray*}
due to the condition $h(x)\ge\ln x$
for all sufficiently large $x$.
This estimate and (\ref{E3.1}) imply
\begin{eqnarray}\label{E3}
E_3 &\le& {\bf E}\{e^{h(\eta)};\eta>g(x)\}
e^{1-\ln g(x)}\nonumber\\
&=& o(1)/g(x) = o(h'(x))
\quad\mbox{ as }x\to\infty,
\end{eqnarray}
by the condition ${\bf E}e^{h(\eta)}<\infty$.
Substituting (\ref{E1}), (\ref{E2}) and (\ref{E3})
into (\ref{decomp.estimate}) and taking into account
the choice (\ref{choice.A}) of $A$, we get
\begin{eqnarray*}
{\mathbf E}e^{h(x+\eta)}
&=& e^{h(x)}{\bf E}e^{h(x+\eta)-h(x)}\\
&\le& e^{h(x)}(1+h'(x)\varepsilon+o(h'(x)))
\quad\mbox{ as }x\to\infty.
\end{eqnarray*}
Since $\varepsilon<0$, the latter estimate
implies ${\bf E}e^{h(x+\eta)} < e^{h(x)}$
for all sufficiently large $x$.
The proof is complete.

{\it Proof} of Lemma \ref{lem.conc2}.
Put $\eta_n=\xi_n-c$. We have ${\bf E}\eta_n<0$ and
${\bf E}e^{h(\eta_n)}<\infty$.
By Lemma \ref{semi.moments.th}, there exists $x_0>0$
such that ${\bf E}e^{h(x+\eta_n)}\le {\bf E}e^{h(x)}$
for $x>x_0$. Then, by monotonicity of $h(x)$ and by
non-negativity of $S_{n-1}$,
\begin{eqnarray*}
{\bf E}e^{h(S_n)}
&\le& {\bf E}e^{h(S_n+x_0)}
= {\bf E}e^{h(S_{n-1}+x_0+c+\eta_n)}
\le {\bf E}e^{h(S_{n-1}+x_0+c)}.
\end{eqnarray*}
Now, by the induction arguments,
${\bf E}e^{h(S_n)}\le e^{h(cn+x_0)}
\le e^{h(cn)}e^{h(x_0)}$.
The proof is complete.

{\bf 4. Proof of Theorem~\ref{main_result}.}
Before starting the proof of Theorem~\ref{main_result},
we formulate the following proposition
from [\ref{DFK2}, Corollary 1]:

\begin{Proposition}\label{general.th}
Let there exist a concave function
$r:{\bf R}^+\to{\bf R}$ such that
${\bf E}e^{r(\xi)}<\infty$ and ${\bf E}\xi e^{r(\xi)}=\infty$.
If $F$ is heavy-tailed and
${\bf E}\tau e^{r(S_{\tau-1})}< \infty$,
then {\rm(\ref{lower.limit.concl})} holds.
\end{Proposition}

We also need two auxiliary technical results.

\begin{Lemma}\label{conc.finite.mom}
Let $\chi\ge 0$ be any random variable.
Then there exists a differentiable concave function
$g:{\bf R}^+\to{\bf R}^+$, $g(0)=0$,
such that $g'(x)\le 1$ for all $x$,
$g(x)\to\infty$ as $x\to\infty$,
and ${\bf E}e^{g(\chi)}<\infty$.
\end{Lemma}

\proof
Consider an increasing sequence $\{x_n\}$
such that $x_0=0$, $x_1=1$, $x_{n+1}-x_n>x_n-x_{n-1}$,
and ${\bf P}\{\chi>x_n\}\le e^{-n}$.
Put $g_1(x_n)=n/2$ and continiously linear
between these points. Then, for any
$x\in(x_n,x_{n+1})$ and $y\in(x_{n+1},x_{n+2})$ we have
$$
g_1'(x)=\frac{1}{2(x_{n+1}-x_n)}
>\frac{1}{2(x_{n+2}-x_{n+1})}=g_1'(y),
$$
so that $g_1$ is concave. By the construction,
$g_1(x)\uparrow\infty$ as $x\to\infty$ and
$g_1'(x)\le1$ where the derivative exists.
Finally,
$$
{\bf E}e^{g_1(\chi)} \le \sum_{n=0}^\infty
e^{g_1(x_{n+1})}{\bf P}\{\chi>x_n\}
\le\sum_{n=0}^\infty e^{(n+1)/2} e^{-n}<\infty.
$$
A procedure of smoothing, say
$g(x)=\int_x^{x+1}g_1(y)dy-\int_0^1g_1(y)dy$,
completes the proof.

\begin{Lemma}\label{from.concave.to.o}
Let $\chi\ge 0$ be a random variable such that,
for some concave function $f:{\bf R}^+\to{\bf R}^+$,
${\bf E}e^{f(\chi)}=\infty$.
Then there exists a concave function
$f_1:{\bf R}^+\to{\bf R}^+$ such that
$f_1\le f$, $f_1(x)=o(x)$ as $x\to\infty$,
and ${\bf E}e^{f_1(\chi)}=\infty$.
\end{Lemma}

\proof
Take $x_1$ so large that
${\bf E}\{e^{\min(\chi,f(\chi))};\chi\le x_1\}\ge1$
and put $f_1(x)=\min(x,f(x))$ for $x\in[0,x_1]$.
Then by induction, for any $n$, we can choose
$x_{n+1}$ such that
$$
{\bf E}\{e^{f_1(x_n)+
\min(n^{-1}f_1'(x_n)(\chi-x_n),f(\chi)-f(x_n))};
\chi\in(x_n,x_{n+1}]\}
\ge 1.
$$
Let $f_1(x)=f_1(x_n)+\min(n^{-1}f_1'(x_n)(x-x_n),f(x)-f(x_n))$
for $x\in(x_n,x_{n+1}]$.
By construction, $f_1$ is concave, $f_1\le f$,
and $f_1'(x_{n+1})\le f_1'(x_n)/n\to 0$ as $n\to\infty$.

{\it Proof} of Theorem \ref{main_result}.
Without loss of generality, assume
that $f(x)\ge\ln x$ for all $x$ and that
$f_2(x)\equiv f(x)-\ln x$ is concave
on the whole posititive half-line.
By Lemma \ref{from.concave.to.o} and
by measure change arguments like in the proof
of Lemma \ref{lm.moments.ext.f1f2}
we may assume from the very beginning that
$$
f(x)=o(x)\quad\mbox{ as }x\to\infty.
$$
Next we state the existence of a concave function
$g:{\bf R}^+\to{\bf R}$ such that
$g(x)\to\infty$ as $x\to\infty$,
$g(x)\le\ln x$ for all sufficiently large $x$,
the difference $\ln x-g(x)$
is a non-decreasing function, and
$$
{\bf E}e^{f(c\tau)+g(c\tau)}<\infty.
$$
Indeed, by Lemma \ref{conc.finite.mom}
and again measure change technique,
there exists a differentiable concave function
$g_1:{\bf R}^+\to{\bf R}^+$
such that $g_1(0)=0$, $g_1(x)\uparrow\infty$, $g_1'(x)\le1$,
and ${\bf E}e^{f(c\tau)+g_1(c\tau)}<\infty$.
Put $g(x)=g_1(\ln(x+1))-1$. Then $g$ is a monotone
function increasing to infinity and
$g(x)\le \ln x$ for all sufficiently large $x$.
In addition,
$$
(\ln x-g(x))'=1/x-g_1'(\ln(x+1))/(x+1) \ge 0,
$$
so that the difference $\ln x-g(x)$
is a non-decreasing function as needed.

Since the function $f_2(x)$ is concave,
by Lemma \ref{lm.moments.ext.f1f2} with
$f_1(x)=\ln x$, there exists a concave function
$h$ such that $h\le f_2$, $h(x)=o(x)$,
${\bf E}\xi e^{h(\xi)}<\infty$ and
${\bf E}\xi e^{h(\xi)+g(\xi)}=\infty$.
Since $\ln x+h(x)+g(x)\le f(x)+g(x)$,
by (\ref{eq.2}) and by the choice of $g$,
\begin{eqnarray}\label{eq.2.ext2}
{\bf E}\tau e^{h(c\tau)+g(c\tau)}<\infty.
\end{eqnarray}
The concave function $r(x)=h(x)+g(x)$ satisfies all
conditions of Proposition~\ref{general.th}.
Indeed, due to the inequality
$g(x)\le\ln x$ for all sufficiently large $x$,
we have ${\bf E}e^{r(\xi)}<\infty$ because
${\bf E}\xi e^{h(\xi)}<\infty$.
It remains to check that
${\bf E}\tau e^{r(S_{\tau-1})}<\infty$.
Since, by (\ref{eq.2.ext2}),
$$
{\bf E}\{\tau e^{r(S_\tau)};S_\tau\le c\tau\}
\le {\bf E}\tau e^{r(c\tau)}<\infty,
$$
it suffices to prove that
\begin{eqnarray*}
{\bf E}\{\tau e^{r(S_\tau)};S_\tau>c\tau\}<\infty.
\end{eqnarray*}
We proceed in the following way:
\begin{eqnarray*}
{\bf E}\{c\tau e^{r(S_\tau)};S_\tau>c\tau\} &=&
\sum_{n=1}^\infty{\bf P}\{\tau=n\}cn{\bf E}\{e^{r(S_n)};S_n>cn\}\\
&=& \sum_{n=1}^\infty{\bf P}\{\tau=n\}e^{g(cn)+\ln(cn)-g(cn)}
{\bf E}\{e^{h(S_n)+g(S_n)};S_n>cn\}.
\end{eqnarray*}
By the monotonicity of the difference $\ln x-g(x)$,
we obtain the following estimate
\begin{eqnarray*}
{\bf E}\{c\tau e^{r(S_\tau)};S_\tau>c\tau\} &\le&
\sum_{n=1}^\infty{\bf P}\{\tau=n\}e^{g(cn)}
{\bf E}\{e^{\ln S_n+h(S_n)};S_n>cn\},
\end{eqnarray*}
Since the function $\ln x+h(x)$ is concave and
$\ln x+h(x)\ge\ln x$, by Lemma \ref{lem.conc2},
$$
{\bf E}e^{\ln S_n+h(S_n)}\le K(c)e^{\ln(nc)+h(cn)}
$$
for some $K(c)<\infty$. Therefore,
\begin{eqnarray*}
{\bf E}\{c\tau e^{r(S_\tau)};S_\tau>c\tau\}
&\le& K(c)\sum_{n=1}^\infty{\bf P}\{\tau=n\}
e^{g(cn)} e^{\ln(cn)+h(nc)}\\
&=& K(c)c{\bf E}\tau e^{h(c\tau)+g(c\tau)} < \infty,
\end{eqnarray*}
from (\ref{eq.2.ext2}).
The proof of Theorem~\ref{main_result} is complete.

{\bf 5. Proof of Theorem \ref{ht.th}.}
Denote by $G$ the distribution function of $c\tau$.

We will construct an increasing concave function
$f:{\bf R}^+\to{\bf R}$ such that
\begin{equation} \label{condi}
{\bf E}\xi e^{f(\xi)}=\infty
\quad \quad \mbox{and}
\quad \quad
{\bf E}\tau e^{f(c\tau)}<\infty.
\end{equation}
Then the desired relation \ref{lower.limit.concl})
will follow by applying Theorem \ref{main_result}.
%in order to obtain the desired relation
%(\ref{lower.limit.concl}).

If $G$ is light-tailed then one can take $f(x)=\lambda x$
for a sufficiently small $\lambda>0$.
>From now on we assume $G$ to be heavy-tailed.

Consider new random variables $\xi_*$
and $\tau_*$ with the following distributions:
\begin{eqnarray*}
{\bf P}\{\xi_*\in dx\}
= \frac{xF(dx)}{{\bf E}\xi}
&\mbox{ and }&
{\bf P}\{\tau_*=n\}
= \frac{n{\bf P}\{\tau=n\}}{{\bf E}\tau}.
\end{eqnarray*}
Denote by $F_*$ and $G_*$ the distributions of $\xi_*$
and $c\tau_*$ respectively.
Then both $F_*$ and $G_*$ are heavy-tailed and
\begin{eqnarray}\label{G.o.F.*}
\overline G_*(x) &=& o(\overline F_*(x))
\quad\mbox{ as }x\to\infty.
\end{eqnarray}
The heavy-tailedness of $G_*$ is equivalent to
the following condition: for any $\varepsilon >0$,
\begin{eqnarray}\label{h.t.G*}
\int_1^\infty\overline G_*(\varepsilon^{-1}\ln x)dx
\equiv \int_0^\infty e^x\overline G_*(x/\varepsilon)dx
&=& \infty.
\end{eqnarray}
%for every $\varepsilon>0$.
In terms of new distributions $F_*$ and $G_*$,
conditions (\ref{condi})
nay be reformulated as follows:
we need to construct an increasing concave function $f$
such that ${\bf E}e^{f(\xi_*)} = \infty$
and ${\bf E}e^{f(c\tau_*)} < \infty$,
or, equivalently,
\begin{equation} \label{condi2}
\int_1^\infty\overline F_*(f^{-1}(\ln x))dx
= \infty
\quad \quad
\mbox{ and }
\quad \quad
\int_1^\infty\overline G_*(f^{-1}(\ln x))dx
< \infty.
\end{equation}
The concavity of $f$ is equivalent to the convexity of its
inverse, $h=f^{-1}$. So, conditions (\ref{condi2}) may be
rewritten as: we have to present an increasing convex function $h$
such that
\begin{equation} \label{condi3}
\int_0^\infty e^x\overline F_*(h(x))dx=\infty
\quad \quad
\mbox{ and }
\quad \quad
\int_0^\infty e^x\overline G_*(h(x))dx<\infty.
\end{equation}

We will construct $h(x)$ as a piece-wise linear
function. For this, we will introduce two increasing
sequences, say $x_n\uparrow\infty$ and
$a_n\uparrow\infty$, and let
\begin{eqnarray*}
h(x) &=& h(x_n)+a_n(x-x_n)
\quad\mbox{ for }x\in(x_n,x_{n+1}].
\end{eqnarray*}
Then the convexity of $f$ will follow from the increase of
$\{a_n\}$.

Put $x_0=0$ and $f(x_0)=0$.
Due to (\ref{G.o.F.*}) and (\ref{h.t.G*}),
we can choose $x_1$ so large that
\begin{eqnarray*}
\frac{\overline F_*(y)}{\overline G_*(y)}
&\ge& 2^1
\end{eqnarray*}
for all $y>x_1$ and
\begin{eqnarray*}
\int_0^{x_1} e^x\overline G_*(h(x_0)+1\cdot(x-x_0))dx
&\ge& 1.
\end{eqnarray*}
Then there exists a sufficiently large
$a_0\ge1$ such that
\begin{eqnarray*}
\int_0^{x_1} e^x\overline G_*(h(x_0)+a_0(x-x_0))dx
&=& 1.
\end{eqnarray*}

Now we use the induction argument to
construct increasing sequences
$\{x_n\}$ and $\{a_n\}$ such that
\begin{eqnarray}\label{frac.F.G.ge}
\frac{\overline F_*(y)}{\overline G_*(y)}
&\ge& 2^{n+1}
\end{eqnarray}
for all $y>x_{n+1}$ and
\begin{eqnarray*}
\int_{x_n}^{x_{n+1}} e^x\overline G_*(h(x))dx
&=& 2^{-n}.
\end{eqnarray*}
For $n=0$ this is already done.
Make the induction hypothesis for some $n\ge1$.
For any $x>x_{n+1}$, denote
\begin{eqnarray*}
\delta(x,a) &\equiv&
\int_{x_{n+1}}^x e^y\overline
G_*(h(x_{n+1}+a(y-x_{n+1})))dy.
\end{eqnarray*}
Due to (\ref{G.o.F.*}) and (\ref{h.t.G*}),
we can choose $x_{n+2}$ so large that
\begin{eqnarray*}
\frac{\overline F_*(y)}{\overline G_*(y)}
&\ge& 2^{n+2}
\end{eqnarray*}
for all $y>x_{n+2}$ and
\begin{eqnarray*}
\delta(x_{n+2},a_n) &\ge& 1.
\end{eqnarray*}
Since the function $\delta(x_{n+2},a)$
continuously decreases to $0$
as $a\uparrow\infty$,
we can choose
$a_{n+1}>a_n$ such that
\begin{eqnarray*}
\delta(x_{n+2},a_{n+1})
&=& 2^{-(n+1)}.
\end{eqnarray*}
Then
\begin{eqnarray*}
\int_{x_{n+1}}^{x_{n+2}}
e^x\overline G_*(h(x))dx
&=& 2^{-(n+1)}.
\end{eqnarray*}
Our induction hypothesis now holds with $n+1$
in place of $n$ as required.

Now the inequalities (\ref{condi3}) follow since,
from the construction of function $h$,
\begin{eqnarray*}
\int_0^\infty e^x\overline G_*(h(x))dx
&=& \sum_{n=0}^\infty
\int_{x_n}^{x_{n+1}}
e^x\overline G_*(h(x))dx\\
&=& \sum_{n=0}^\infty 2^{-n}<\infty.
\end{eqnarray*}
and, by (\ref{frac.F.G.ge}),
\begin{eqnarray*}
\int_0^\infty e^x\overline F_*(h(x))dx
&=& \sum_{n=0}^\infty
\int_{x_n}^{x_{n+1}}
e^x\overline F_*(h(x))dx\\
&\ge& \sum_{n=0}^\infty
2^n\int_{x_n}^{x_{n+1}}
e^x\overline G_*(h(x))dx\\
&=& \sum_{n=0}^\infty 2^n 2^{-n}=\infty.
\end{eqnarray*}
The proof of Theorem \ref{ht.th} is complete.

{\bf 6. Proof of Theorem \ref{lt.th}.}
We apply the exponential change of measure with
parameter $\widehat\gamma$ and consider the distribution
$G(du)=e^{\widehat\gamma u}F(du)/\varphi(\widehat\gamma)$
and the stopping time $\nu$ with the distribution
${\bf P}\{\nu=k\}=\varphi^k(\widehat\gamma){\bf P}\{\tau=k\}/
{\bf E}\varphi^\tau(\widehat\gamma)$.
Then it was proved in [\ref{DFK2}, Lemma 3] that
\begin{eqnarray}\label{3.3}
\liminf_{x\to\infty}
\frac{\overline{G^{*\nu}}(x)}{\overline G(x)}
&\ge& \frac{1}{{\bf E}\varphi^{\tau-1}(\gamma)}
\liminf_{x\to\infty}
\frac{\overline{F^{*\tau}}(x)}{\overline F(x)}.
\end{eqnarray}
>From the definition of $\widehat\gamma$,
the distribution $G$ is heavy-tailed.
Let us prove that
\begin{eqnarray}\label{c.nu.o}
{\bf P}\{c\nu>x\} &=& o(\overline G(x))
\quad\mbox{ as }x\to\infty.
\end{eqnarray}
Indeed, put
$\lambda\equiv\ln\varphi(\widehat\gamma)>0$;
then
\begin{eqnarray}\label{c.nu.rep}
{\bf P}\{c\nu>x\} &=&
\frac{1}{{\bf E}\varphi^\tau(\widehat\gamma)}
\sum_{k>x/c}e^{\lambda k}{\bf P}\{\tau=k\}\nonumber\\
&\le& \frac{1}{{\bf E}\varphi^\tau(\widehat\gamma)}
\int_{x/c}^\infty e^{\lambda y}{\bf P}\{\tau\in dy\}.
\end{eqnarray}
Integration by parts implies
\begin{eqnarray*}
\int_{x/c}^\infty e^{\lambda y}{\bf P}\{\tau\in dy\}
&=& -e^{\lambda y}{\bf P}\{\tau>y\}\Big|_{x/c}^\infty+
\lambda\int_{x/c}^\infty e^{\lambda y}{\bf P}\{\tau>y\}dy\\
&=& e^{\lambda x/c}{\bf P}\{c\tau>x\}+
\frac{\lambda}{c}\int_x^\infty e^{\lambda y/c}
{\bf P}\{c\tau>y\}dy,
\end{eqnarray*}
because ${\bf E}\varphi^\tau(\widehat\gamma)<\infty$ and, thus,
$e^{\lambda y}{\bf P}\{\tau>y\}\to0$ as $y\to\infty$.
Now applying the condition (\ref{G.=o.F})
we obtain that the latter sum is of order
\begin{eqnarray*}
o\Bigl(e^{\lambda x/c}\overline F(x)+
\frac{\lambda}{c}\int_x^\infty e^{\lambda y/c}\overline F(y)dy
\Bigr)
&=& o\Bigl(\int_x^\infty e^{\lambda y/c}F(dy)\Bigr)
\quad\mbox{ as }x\to\infty.
\end{eqnarray*}
Together with (\ref{c.nu.rep}) it implies (\ref{c.nu.o}).
Therefore, by Theorem \ref{ht.th} we have the equality
\begin{eqnarray*}
\liminf_{x\to\infty}
\frac{\overline{G^{*\nu}}(x)}{\overline G(x)}
&=& {\bf E}\nu
=\frac{{\bf E}\tau\varphi^\tau(\widehat\gamma)}
{{\bf E}\varphi^\tau(\widehat\gamma)},
\end{eqnarray*}
and, due to (\ref{3.3}),
\begin{eqnarray}\label{liminf.le}
\liminf_{x\to\infty}
\frac{\overline{F^{*\tau}}(x)}{\overline F(x)}
&\le& {\bf E}\tau\varphi^{\tau-1}(\widehat\gamma).
\end{eqnarray}

The result now follows from Lemma .

\section*{\normalsize References}

\newcounter{bibcoun}
\begin{list}{\arabic{bibcoun}.}{\usecounter{bibcoun}\itemsep=0pt}
\small

\item\label{CNW}
Chover, J., Ney, P. and Wainger, S., 1973.
Functions of probability measures.
J. Anal. Math. 26, 255--302.

\item\label{Cline}
Cline, D., 1987.
Convolutions of distributions with
exponential and subexponential tailes.
J. Aust. Math. Soc. 43, 347--365.

\item\label{DFK2}
Denisov, D., Foss, S. and Korshunov, D.
On lower limits and equivalences for tails of randomly stopped sums.
To appear in Bernoulli.

\item\label{EG}
Embrechts, P. and Goldie, C. M., 1982.
On convolution tails.
Stochastic Process. Appl. 13, 263--278.

\item\label{FK}
Foss, S. and Korshunov, D., 2007.
Lower limits and equivalences for convolution tails.
Ann. Probab. 35, 366-–383.

\item\label{FS}
Foss, S. and Sapozhnikov, A., 2004.
On the existence of moments for the busy period
in the single-server queue.
Math. Oper. Research 29, 592--601.

\item\label{Pakes}
Pakes, A. G., 2004.
Convolution equivalence and infinite divisibility.
J. Appl. Probab. 41, 407--424.

\item\label{Rogozin}
Rogozin, B. A., 2000.
On the constant in the definition of
subexponential distributions.
Theory Probab. Appl. 44, 409--412.

\item\label{Rudin}
Rudin, W., 1973.
Limits of ratios of tails of measures.
Ann. Probab. 1, 982--994.
\end{list}

\end{document}